\documentclass[10pt]{article}

\usepackage{amssymb}
\usepackage{amsthm}
\usepackage{amsmath}
\usepackage{amscd}

\allowdisplaybreaks

\newcommand{\R}{\mathbb R}

\newcommand{\eps}{\varepsilon}

\newcommand{\dif}{\mathrm{d}}
\newcommand{\Dif}{\mathrm{D}}

\newcommand{\approxsol}{\tilde S_{t}}
\newcommand{\cobg}{C^{0,\alpha}_{\gamma - 2}}
\newcommand{\ctbg}{C^{2,\alpha}_\gamma}
\newcommand{\ckbg}{C^{k,\alpha}_\gamma}

\newcommand{\Btbg}{\mathcal{B}^{2,\alpha}_{\gamma} (\tilde S_{t})}
\newcommand{\Bobg}{\mathcal{B}^{0,\alpha}_{\gamma - 2} (\tilde S_{t})}
\newcommand{\spqa}{C_t}

\newcommand{\dopq}{O(p+1) \times O(q+1)}

\newtheorem{thm}{Theorem}

\newtheorem{cor}[thm]{Corollary}
\newtheorem{prop}[thm]{Proposition}
\newtheorem*{nonumthm}{Theorem}

\theoremstyle{definition}
\newtheorem{defn}[thm]{Definition}

\newcommand{\mylabel}
    {\label}

\begin{document}

\title{Generalized Doubling Constructions\\  for Constant Mean Curvature Hypersurfaces in $S^{n+1}$}

\author{
\begin{minipage}{3.125in}
    \begin{center}
        Adrian Butscher \\ University of Toronto at Scarborough \\ email: \ttfamily butscher@utsc.utoronto.ca
    \end{center}
\end{minipage}
\begin{minipage}{3.125in}
\begin{center}
        Frank Pacard \\ Universit\'e de Paris XII \\ email: \ttfamily pacard@univ-paris12.fr
    \end{center}
\end{minipage}\\
\rule{1ex}{0ex}
}
\maketitle

\begin{abstract}
The sphere $S^{n+1}$ contains a simple family of constant mean curvature
(CMC) hypersurfaces of the form $C_t := S^p ( \cos t ) \times S^q ( \sin t  )$ for $p+q = n$ and $t \in ( 0 , \tfrac{\pi}{2})$ called the generalized Clifford hypersurfaces. This paper demonstrates that new, topologically non-trivial CMC hypersurfaces resembling a pair of neighbouring generalized Clifford tori connected to each other by small catenoidal bridges at a sufficiently symmetric configuration of points can be constructed by perturbative PDE methods. That is, one can create an approximate solution by gluing a rescaled catenoid into the neighbourhood of each point; and then one can show that a perturbation of this approximate hypersurface exists which satisfies the CMC condition. The results of this paper generalize those of the authors in \cite{mepacard1}.
\end{abstract}

\renewcommand{\baselinestretch}{1.25}
\normalsize

\section{Introduction and Statement of Results}

\paragraph*{CMC hypersurfaces.} A constant mean curvature (CMC) hypersurface $\Sigma$ contained in an ambient Riemannian manifold $X$ of dimension $n+1$ has the property that its mean curvature with respect to the induced metric is constant.  This property ensures that $n$-dimensional area of $\Sigma$ is a critical value of the area functional for hypersurfaces of $X$ subject to an enclosed-volume constraint.  Constant mean curvature hypersurfaces have been objects of great interest since the beginnings of modern differential geometry.  Classical examples of non-trivial CMC surfaces in three-dimensional Euclidean space $\R^3$ are the sphere, the cylinder and the Delaunay surfaces, and for a long while these were the only known CMC surfaces.  In fact, a result of Alexandrov \cite{alexandrov} states that the only compact, connected, embedded CMC surfaces in $\R^3$ are spheres.

In recent decades, the theory of CMC surfaces in $\R^3$ has progressed considerably. However, the corresponding picture amongst CMC hypersurfaces of higher dimension or in other ambient manifolds is not nearly as rich, due in part to the absence of the Weierstra\ss -type representation or the Lawson associated surface construction that are available in $\R^3$.  There is a certain amount of literature on CMC hypersurfaces in hyperbolic space \cite{bryanthyper,pacardpimentel,saerptoubiana,umeharayamada}; but due to the non-compactness of hyperbolic space, this theory can be considered not such a vast departure from the theory of CMC hypersurfaces in $\R^{n+1}$.  Much less is known when the ambient space is the sphere.  The classically known examples in $S^{n+1}$ are the hyperspheres obtained from intersecting $S^{n+1}$ with hyperplanes, and the so-called generalized Clifford tori which are products of lower-dimensional spheres of the form 
\[
C_t : = S^p \left( \cos t \right) \times S^q \left( \sin t \right)
\] 
for $p+q = n$ and $t \in (0 , \tfrac{\pi}{2})$.  These are embedded hypersurfaces in $S^{n+1}$ with constant mean curvature equal to $H_t : =  q \, \cot t - p \tan t$.  There are few other examples, and no general methods for the construction of CMC surfaces in $S^{n+1}$.  However, the method of \emph{gluing}, in which a CMC hypersurface is constructed by pasting together simple building blocks, is a successful technique in the $\R^3$ setting and can be attempted in $S^{n+1}$.  This is because many of the operations involved in a gluing construction --- such as forming connected sums using small bridging surfaces near a point of mutual tangency --- are all \emph{local} and thus have straightforward generalizations to other ambient manifolds.   

When $n=2$ and hence $p=q=1$, Butscher and Pacard have proven in \cite{mepacard1} that in $S^3$, it is possible to construct new examples of embedded, higher-genus CMC surfaces of $S^3$, with small but non-zero mean curvature, by \emph{doubling} the unique minimal Clifford torus $C_{\tfrac{\pi}{4}}$ in the family of Clifford tori of $S^3$.  That is, these new surfaces are small perturbations of two parallel translates of $C_{\tfrac{\pi}{4}}$ which are glued together at a sub-lattice of points by means of small catenoidal bridging surfaces.  The two parallel translates are a distance $\eps$ apart and the mean curvature of the doubled surfaces is given by $H = \tfrac{1}{2} \, \cot (\tfrac{\pi}{4}+ \eps)$.  When $\eps$ tends to zero, the doubled surface converges, away from the points where the catenoids are glued, to two copies of  $C_{\tfrac{\pi}{4}}$.  These surfaces are in a certain sense compact analogues of the doubly periodic CMC surfaces in $\R^3$ constructed by Ritor\'e \cite{ritore} and Gro\ss{e}-Brauckmann and Karcher \cite{grossebrauckmann,karcher}.

\paragraph*{The generalized doubling construction.}  This paper generalizes the Butscher-Pacard construction to the sphere $S^n$.  The family of generalized Clifford hypersurfaces $C_t$ is also a foliation a tubular neighbourhood of the minimal hypersurface $C_{t_*}$, with 
\[
\tan t_* : = \sqrt{\tfrac{q}{p}} ,
\]
having parallel leaves.  Thus two parallel translates of $C_{t_*}$, located on each side, can be connected together at a symmetric configuration of points, called the \emph{gluing points}, by means of $n$-dimensional catenoidal bridges.  The resulting hypersurface, henceforth called $\approxsol$, can be constructed with various kinds of non-trivial topology, depending on the number of gluing points.  Once again, $\approxsol$ is only approximately CMC and must be perturbed to achieve constant mean curvature.  This perturbation is in general obstructed due to the existence of non-trivial elements of the kernel of the linearized mean curvature operator of the constituents of $\approxsol$, called \emph{Jacobi fields}, whose effect is to prevent the linearized CMC equation from being bijective with bounded inverse.  As in the $S^3$ case, the way to avoid the obstructions is to impose additional symmetries on the approximate solution that are not possessed by the Jacobi fields.  That is, if the gluing points are chosen with sufficient symmetry and $\approxsol$ is perturbed in a way which respects these symmetries, then one can show that the Jacobi fields are absent and the CMC equation can be controllably inverted.

The most economical way of encoding the symmetries necessary for the construction outlined above is {\it via} finite subgroups of symmetries of $S^{n+1}$ that preserve $C_{t_*}$.   The use of symmetry groups generalizes the sub-lattice of the torus $S^1 \times S^1$ used in the $S^3$ case.  The $S^{n+1}$ case requires a more sophisticated choice in part because at least one of the spherical factors of $S^p \times S^q$ is itself of higher dimension, in which case the analogue of a sub-lattice is not natural. The symmetry condition necessary for the proof of the present theorem can be explained as follows.  Note first that the full group of symmetries preserving $C_t$ is exactly $\dopq$ acting diagonally on $\R^{n+2} = \R^{p+1} \times \R^{q+1}$.  The finite subgroups that we are interested  will be of the form 
\[
G \subset \{ (\sigma_{p+1} , \sigma_{q+1}) \: : \: \sigma_s \in O(s) \: \mbox{for $s = p+1, q+1$} \}.
\]  
We will assume that $G$ contains the element  $\rho := (\rho_{p+1}, \rho_{q+1})$ where $\rho_N \in O(N+1)$ is the reflection symmetry across the $x_1 =0$ axis, namely 
\[
\rho_N ((x_1, x_2,\ldots, x_{N+1}))= ((x_1, - x_2, \ldots, - x_{N+1})).
\] 
Next, we define the point $\mu_0 \in C_{t_*}$ to be \[
\mu_0 := (\sqrt{\tfrac{p}{n}} , 0, \ldots, 0, \sqrt{\tfrac{q}{n}}, \ldots, 0)
\] 
and the set $\Lambda \subset C_{t_*}$ to be the orbit of $\mu_0$ under $G$. We denote the cardinality of $\Lambda$ by $m_\Lambda$. 

\begin{thm}
\label{th:1}
Assume that there are no numbers $a^{kl} \in \R$ (not all equal to $0$) such that the function 
\[
(x,y) \in \R^{p+1} \times \R^{q+1} \longmapsto \sum_{k,l} a^{kl} \, x_k \, y_l
\]
is $G$-invarariant. Then for all $t$ close enough to $t_*$, there exists a smooth, embedded, CMC hypersurface $S_t$ with the following properties.
    \begin{enumerate}
\item The hypersurface $S_t$ is invariant under the action of $G$.
\item The hypersurface $S_t$ is topologically equal to the connected sum of two copies of $(S^p \times S^q)$ at $m_\Lambda$ points.
\item The mean curvature of $S_t$ is equal to $H_t : =  q \, \cot t - p \, \tan t$.
        \item Away from a neighbourhood of $\Lambda$, the hypersurface $S_t$ is a perturbation of two hypersurfaces in the family $C_t$ located on either side of $C_{t_*}$.

        \item In a neighbourhood of each point in $\Lambda$, the hypersurface $S_t$ is a perturbation of a truncated, rescaled $n$-dimensional catenoid whose ends are attached to the hypersurfaces described in (4).

        \item  As $t$ tends to $t_*$ then $S_t$ converges in $C^\infty$ topology to two copies of $C_{t_*}$ away from the points of $\Lambda$.
    \end{enumerate}
\end{thm}

The proof of this theorem proceeds in a parallel fashion to the proof of the version valid in $S^3$ that is given in \cite{mepacard1}.  This begins by generalizing the initial construction of \cite{mepacard1}, whereby two normal translates of $C_{t_*}$ separated by a small amount are first glued together at the gluing points using small necks of height equal to the distance between the translates.  This is the \emph{approximate solution}  $\approxsol$.    It makes sense that such a construction is possible only if $C_{t_*}$ is the unique minimal generalized Clifford hypersurface since then normal translates to either side of $C_{t_*}$ can be chosen which have opposite mean curvature.  It will then be shown that $\approxsol$ can be perturbed to have exactly constant mean curvature.  

\section{Examples}

In this section we given examples of the application of Theorem~\ref{th:1}. Basically, we give examples of groups for which the main assumption of Theorem~\ref{th:1} is fulfilled. In what follows, it is easiest to describe the elements of $G$ through their action on ${\mathbb R}^{n+2}$.

\paragraph*{Example 1.} In the lowest-dimensional case $n=2$ and $p=q=1$,  considered by Butscher and Pacard in \cite{mepacard1}, two normal translates of $C_{\tfrac{\pi}{4}}$ are glued together at a sub-lattice of points.  We choose $\tau_j : = (\alpha_j, \beta_j) \in {\mathbb R}^2$ with $j=1,2$ so that the lattice ${\mathbb Z} \,  \tau_1 + {\mathbb Z} \, \tau_2$ contains the lattice $2 \, \pi \,  {\mathbb Z}^2$.  Let $R_\alpha $ denote the rotation of angle $\alpha$ in ${\mathbb R}^2$. We consider the group $G$ generated by  the elements of $O(4)$ whose actions are given by
\begin{align*}
	\sigma_j  (x_1, x_2 ,y_1, y_2 ) &:=  \left(  R_{\alpha_j}  (x_1, x_2) ,  R_{\beta_j} (y_1, y_2) \right) \\
	\intertext{for $j=1,2$ and}
	\rho  (x_1, x_2 ,y_1, y_2 )  &:=  (x_1,  - x_2 ,y_1, - y_2 ) \, .
\end{align*}
It is proven in \cite{mepacard1} that the condition on $\tau_j$ that ensures that  is the following: the lattice ${\mathbb Z}  \, \tau_1 + {\mathbb Z} \, \tau_2$ is not contained in $\{(\alpha, \beta) \in {\mathbb R}^2 \: : \: \alpha - \beta \equiv 0 \: [2\pi]\}$ or in $\{(\alpha, \beta) \in {\mathbb R}^2 \: : \: \alpha + \beta \equiv 0 \: [2\pi]\}$.

\paragraph*{Example 2.} The previous example extends {\it verbatim} to any dimension. We consider the group $G$ generated by  the elements of $O(n+2)$ whose actions are given by
\begin{align*}
	\sigma_j  (x_1, \ldots , x_{p+1},y_1, \ldots, y_{q+1}) &:=  \left(  R_{\alpha_j}  (x_1, x_2) ,  x_3, \ldots, x_{p+1},  R_{\beta_j} (y_1, y_2), y_3, \ldots , y_{q+1}\right)  \\
	\intertext{for $j=1,2$ and}
	\rho (x_1, \ldots , x_{p+1},y_1, \ldots, y_{q+1}) &:=  \left( \rho_{p+1} (x), \rho_{q+1}(y) \right)  
\end{align*}
as well as by the $2^{(p-1)(q-1)}$ elements whose action is given by
\[
\rho_{\pm, \ldots, \pm} (x_1, \ldots , x_{p+1},y_1, \ldots, y_{q+1}): =  \left( x_1, x_2, \pm x_3 , \ldots , \pm x_{p+1} , y_1, y_2, \pm y_3 , \ldots , \pm y_{q+1}\right) \, .
\]
In this case, the only function $(x,y) \in \R^{p+1} \times \R^{q+1} \longmapsto \sum_{k,l} a^{kl} \, x_k \, y_l$ that is invariant under the action of $\rho_{\pm, \ldots , \pm}$ is of the form $f(x,y)= \sum_{k,l \in \{1,2\}} a^{kl} \, x_k \, y_l$ and checking that this function is identically equal to $0$ reduces to what is done in \cite{mepacard1}.

\paragraph*{Example 3.}  A third important class of examples is the one where the group $G$ contains the $2^{p\, q}$ elements of $O(n+2)$ whose actions are given by
\[
\tilde \rho_{\pm, \ldots, \pm} (x_1, \ldots , x_{p+1},y_1, \ldots, y_{q+1}) := \left( x_1, \pm x_2, \pm x_3 , \ldots , \pm x_{p+1} , y_1, \pm y_2, \pm y_3 , \ldots , \pm y_{q+1}\right) 
\]
and the orbit of $\mu_0$ by $G$ is not included in $\{\pm \mu_0\}$.
In this case the only function $(x,y) \in \R^{p+1} \times \R^{q+1} \longmapsto \sum_{k,l} a^{kl} \, x_k \, y_l$ that is invariant under the action of the $\tilde \rho_{\pm, \ldots \pm}$ is of the form $f (x,y)= a^{11} \, x_1 \, y_1$ 
and hence has to be identically equal to $0$ if the orbit of $\mu_0$ by $G$ contains more than $\pm \mu_0$.  This is because the value of $f$ at $\mu_0$ is equal to $a^{11}\, \frac{\sqrt{p\, q}}{n}$ on $C_t$, which is a value that is only achieved at $\pm \mu_0$ if $a^{11}\neq 0$.

\section{The Building Blocks of the Doubling Construction}

The purpose of this section is to carefully describe of the building blocks that will be assembled to construct the approximate solution $\approxsol$ --- the generalized Clifford hypersurfaces in $S^{n+1}$ and the generalized catenoid in $\R^{n+1}$.  Since the proof of the Theorem~\ref{th:1} hinges on being able to rule out the existence of Jacobi fields on these building blocks, careful attention will be paid to understanding the Jacobi fields in each case.  Begin with the following characterization of the origin of the Jacobi fields.

\subsection{The Mean Curvature Operator and its Jacobi Fields}

Let $\Sigma$ be a closed hypersurface in a Riemannian manifold $X$ with mean curvature $H_\Sigma$, second fundamental form $B_\Sigma$ and unit normal vector field $N_\Sigma$. The linearization of the mean curvature operator on the space of normal graphs over $\Sigma$ is given by 
\[
{\mathcal L}_\Sigma : = \Dif H_\Sigma (0)  = \Delta_\Sigma  + \Vert B_\Sigma \Vert^2 + {\mathrm{Ric}}(N_\Sigma, N_\Sigma)  
\] 
where $\Delta_\Sigma$ is the Laplace operator of $\Sigma$ and ${\mathrm{Ric}}$ is the Ricci tensor of $X$.  Recall that if $R_t$ is a one-parameter family of isometries of $X$ with deformation vector field $V = \left. \frac{\dif }{\dif t} \right|_{t=0} R_t$, then one obtains a Jacobi field because the function $\langle V, N \rangle$ is in the kernel of $\mathcal L_\Sigma$. When $\Sigma$ is a hypersurface in the ambient space $X=S^{n+1}$, the linearized mean curvature reads 
\[
{\mathcal L}_\Sigma  = \Delta_\Sigma + \Vert B_\Sigma \Vert^2 + n   
\] 
and the isometries of $S^{n+1}$ are simply the $SO(n+2)$-rotations of the ambient $\R^{n+2}$.  Thus there is at most an $(n+2)(n+1)/2$-dimensional space of such `geometric' Jacobi fields of $\Sigma$.  

\subsection{Generalized Clifford Hypersurfaces in \boldmath{$S^{n+1}$}}

\paragraph*{Definition and basic properties.} 

    Let $p$, $q$ and $n \geq 3$ be fixed positive integers such that $p + q  = n$.   The \emph{generalized Clifford hypersurfaces} in $S^{n+1}$ are defined by 
\[
C_t :=  \left\{ (x,y) \in \R^{p+1} \times \R^{q+1} \: : \: \Vert x \Vert = \cos t \; \mbox{ and } \; \Vert y \Vert = \sin t \right\}
\]
for any $t \in (0,\tfrac{\pi}{2})$.
    
Each generalized Clifford hypersurface $C_t$ is topologically equivalent to the product $S^p \times S^q$ and is embedded in $S^{n+1}$.  The following assertions about the geometry of $C_t$ are easy to verify.  First, the induced metric of $\spqa$ is given by
\[
g_t = \cos^2 t \,  g_{S^p} + \sin^2 t \, g_{S^q}
\]
where $g_{S^p}$ and $g_{S^q}$ are the standard metrics on the unit spheres $S^p \subseteq \R^{p+1}$ and $S^q \subseteq \R^{q+1}$, respectively.  The unit normal vector field of $C_t$ is chosen to be 
\[
N_t : =   \sin t \, P_{x} -  \cos t \,  P_{y}
\]
where $P_{x}$ and $P_{y} $ are the position vector fields of $\R^{p+1} \times \{0\} \subset {\R}^{n+2}$ and $\{0\}Ê\times \R^{q+1} \subset {\R}^{n+2}$ respectively.  The second fundamental form  of $C_t$ is given by 
\[
B_t : =  \cos t \, \sin t  \, \left( g_{S^q} -g_{S^p} \right).
\]
Observe that the mean curvature is equal to  
\[
H_t :=   q \,  \cot t  - p \, \tan t \, .
\]
In particular, if $t_* \in (0, \tfrac{\pi}{2})$ is defined by $\tan^2 t_* = \tfrac{q}{p}$, then $C_{t_*}$ has zero mean curvature.  Finally, the linearized mean curvature  operator of $C_t$ is given by 
\[
\mathcal{L}_{t} : = \tfrac{1}{\cos^2 t} \,  (\Delta_{S^p} +p) + \tfrac{1}{\sin^2 t} \, (\Delta_{S^q} +  q) 
\]
where $\Delta_{S^p}$ and $\Delta_{S^q}$ are the Laplacians of $g_{S^p}$ and $g_{S^q}$, respectively.

\paragraph*{Analytic properties of the Jacobi operator.}  The following proposition gathers the necessary information about the Jacobi fields of $C_{t_*}$.

\begin{prop}
    \mylabel{prop:clifjacobi}
    The non-trivial Jacobi fields of $C_{t_*}$ are generated by the $pq$-dimensional subgroup of rotations of $\R^{n+2}$ breaking the $\R^{p+1} \times \R^{q+1}$ splitting. They are  the restriction to $C_{t_*}$ of functions  of the form 
\[
(x,y) \in \R^{p+1} \times \R^{q+1} \longmapsto \sum_{k,l} a^{kl} \, x_k \, y_l \in \R.
\]
where $a^{kl} \in {\mathbb R}$.
\end{prop}

\begin{proof}
Recall that the eigenvalues of $\Delta_{S^N}$ are given by $- j\, (N-1+j)$ for $j \in {\mathbb N}$.  We denote by $E_N$ the eigenfunctions of $\Delta_{S^N}$ associated to the eigenvalue $-N$. Recall that the $E_N$ are the restriction to $S^N$ of linear functions. 

Therefore, the eigenvalues of ${\mathcal L}_{t_*}$ are given by
\[
\lambda_{ij} = - \tfrac{1}{\cos^2 t_*}\, (i^2+ i\, (p-1) - p) - \tfrac{1}{\sin^2 t_*}\, (j^2+ j \, (q-1) - q)
\]
for $i,j \in {\mathbb N}$. Obviously, $\lambda_{ij} < 0$ when $i,j \geq 1$ and $(i,j)\neq (1,1)$. Also $\lambda_{11}=0$ and the corresponding eigenspace is spanned by functions of the form $E_p \, E_q$. Finally 
\[
\lambda_{i0}  = - \tfrac{1}{\cos^2 t_*}\, (i^2+ i\, (p-1) - p) + \tfrac{1}{\sin^2 t_*}\, q = - \tfrac{1}{\sin^2 t_*}\, (i^2+ i\, (p-1) - 2p) 
\]
since $\frac{p}{\cos^2 t_*} =  \frac{q}{\sin^2 t_*}$ and this quantity is never $0$. Similarly $\lambda_{0j} \neq 0$ for all $j\in {\mathbb N}$.
\end{proof}

\noindent Most important for our purposes is the following simple consequence.

\begin{cor}
\label{cor:nojacobi}
Under the assumption of Theorem~\ref{th:1}, there are no Jacobi field on $C_{t_*}$ that is invariant under the action of $G$.
\end{cor}

\subsection{The Generalized Catenoid in \boldmath{$\R^{n+1}$}}

\paragraph*{Definition and basic properties.} In the lowest-dimensional case considered in \cite{mepacard1}, the necks used to glue together two neighbouring Clifford tori were truncations of the standard catenoid in $\R^3$, re-scaled to a small size, and embedded in $S^3$ at the gluing points using canonical coordinate charts.  The appropriate neck in the present higher-dimensional case should then just be the higher-dimensional analogue of the standard catenoid, namely the unique, cylindrically symmetric, minimal hypersurface in $\R^{n+1}$.

The generalized catenoid in $\R^{n+1}$ is the hypersurface $K$ parameterized by
\[
( s , \Theta) \in {\mathbb R}Ê\times S^{n-1} \longmapsto ( \phi(s) \,  \Theta , \psi (s) ) \in {\R}^{n+1}
\]
where  
\begin{equation}
\label{eq:phi}
\phi(s) : =  (\cosh (n-1) s)^{\frac{1}{n-1}} 
\qquad \mbox{and} \qquad \psi (s) : =  \int_0^s \, \phi^{2-n} (t) \, dt \, .
\end{equation}
The geometric features of the generalized catenoid that will be relevant later on are as follows.  The induced metric of $K$ is 
\[
g_K : = \phi^2 \, (\dif s^2 + g_{S^{n-1}}) \, .
\]
The unit normal vector field of $K$ is chosen to be
\[
N_K : = - \phi^{1-n} \, P_\Theta + \partial_s \log \phi \, \partial_{x_{n+2}} \, ,
\]
where $P_\Theta$ is the position vector field in ${\mathbb R}^n \times \{0\} \subset  {\mathbb R}^{n+1}$ evaluated at the point $\Theta \in S^{n-1}$. Then the second fundamental form of $K$ is given by 
\[
B_K : = \phi^{2-n} \, ((1-n) \, \dif s^2 + g_{S^{n-1}})
\] 
and its mean curvature vanishes.  Finally, the Jacobi operator of $K$ is given explicitely by
\[
\mathcal{L}_K  :=  \phi^{-n} \, {\partial_s} \left( \phi^{n-2}  {\partial_s} \right) + \phi^{-2}Ê\, \Delta_{S^{n-1}} + n(n-1) \, \phi^{-2n} \, .
\]
 
\paragraph*{Analytic properties of the Jacobi operator.}  Analytic obstructions for inverting the mean curvature operator on a hypersurface consisting of several large pieces connected by small necks also arise from the non-trivial Jacobi fields of supported in the neck regions.  Thus it is just as important to understand the Jacobi fields on the generalized catenoid in greater detail.

\begin{prop}
    \mylabel{prop:catjacobi}
    Assume that $\delta < 0$ is fixed. Then there is no nontrivial Jacobi field of $K$ that is bounded by a constant times $(\cosh s)^\delta$ and is invariant under the action of the symmetry $(s, \Theta )\longmapsto (s, - \Theta)$.
\end{prop}

\begin{proof}
We consider the eigenfunction decomposition of any Jacobi field
\[
f(s, \cdot )  = \sum_{j=0}^\infty f_j (s, \cdot )
\]
where $\Delta_{S^{n-1}} f_j(s, \cdot) = - j(n-2+j) \, f_j(s, \cdot)$. Since we have assumed that $f(s, -\Theta) = f(s, \Theta)$ many components are equal to $0$ and in particular, we have $f_1 \equiv 0$.  

It follows from \cite{Kab-Pac} that  $f_j \equiv 0$, for all $j \geq 2$ since we have assumed that $|f| \leq C \, (\cosh s)^\delta$ for some  $\delta  < 2$.  Let us briefly remind the reader how this result is proven. First, the geometric Jacobi fields associated to horizontal translations show that the function $u = \phi^{1-n}$ is a solution of 
\[
\phi^{-n} \, {\partial_s} \left( \phi^{n-2}  {\partial_s}  u \right) - \phi^{-2}Ê\,  u + n(n-1) \, \phi^{-2n} u=0 \, .
\]
Now consider $f_j(s, \theta)$ that we decompose on a basis of the $j^{\mathrm{th}}$ eigenspace of $\Delta_{S^{n-1}}$.  The coefficients $u_j^{(\ell)}$ of this decomposition only depend on $s$ and are solutions of 
\[
\phi^{-n} \, {\partial_s} \left( \phi^{n-2}  {\partial_s}  u_j^{(\ell)} \right) - j(n-2+j) \, \phi^{-2}Ê\,  u_j^{(\ell)} + n(n-1) \, \phi^{-2n} u_j^{(\ell)} =0 \, .
\]
Inspection of the possible behaviours of $u_j^{(\ell)}$ shows that, since $j \geq 2$ and $u_j^{(\ell)}$ is bounded by a constant times $(\cosh s)^\delta$ for some $\delta < 2$, then $u_j^{(\ell)}$ is bounded by a constant times  $(\cosh s)^{2-n-j}$. Then, the function $u$, which does not change sign and decays like $(\cosh s)^{1-n}$ at $\pm \infty$, can be used as a barrier to prove that $u_j^{(\ell)} \equiv 0$.
 
The function $f_0$ does not depend on $\theta$ and hence is a solution of some homogeneous second order ordinary differential equation. Two independent solutions of the equations are known since they correspond to geometric Jacobi fields associated to vertical translation and dilation. These are explicitly given by
\[
f^{(1)}_0(s) : = \partial_s \log \phi \qquad \mbox{and} \qquad f^{(2)}_0(s) : = \psi  \, \partial_s \log \phi - \phi^{2-n}
\]
and one checks that no linear combination of these two functions decays exponentially at both $\pm \infty$. This completes the proof of the result.
\end{proof}

\section{The Approximate Solution}

The previous section of this paper described the building blocks of the gluing construction that will be deformed into a CMC hypersurface of $S^{n+1}$.  This section shows in technical detail how these building blocks will be assembled.  We keep the notations of the introduction.  The gluing construction will consist of two generalized Clifford hypersurfaces lying at a small distance on either side of $C_{t_*}$ and glued together at the admissible collection of points $\Lambda$.  The actual gluing will be made using truncated, re-scaled, generalized catenoids embedded into a neighbourhood of each $\muÊ\in \Lambda$ by means of canonical parameterization for a neighbourhood of each point $\mu \in \Lambda$.  We begin by describing this parametrization.

\subsection{Adapted Local Coordinates for \boldmath{$S^{n+1}$}}

We first introduce \emph{toroidal coordinates} for a tubular neighbourhood of $C_{t_*}$.  The coordinate embedding of these coordinates is defined \emph{via} the inverse of the parameterization $\Xi : S^p \times S^q \times  (0, \frac{\pi}{2}) \rightarrow S^{n+1}$ given by
\begin{equation}
    \mylabel{eqn:spherecoords}
    \Xi (z , v ) = \left( \cos v  \, \Theta^{(p)} , \sin  v  \,  \Theta^{(q)} \right)
\end{equation}
for $v \in (0, \frac{\pi}{2})$ and $z := (\Theta^{(p)}, \Theta^{(q)}) \in S^p \times S^q$.  Thus $\Xi$ parameterizes a neighbourhood of $C_{t_*}$ in  $S^{n+1}$.  The local geometry of $S^n$ near $C_{t_*}$ can be completely expressed in the toroidal coordinates.  For instance, the metric is given by
\begin{equation}
    \mylabel{eqn:localmetric}
    \Xi^\ast g_{S^{n+1}} = \dif v^2 + \cos^2 v \,  g_{S^p} + \sin^2 v \,  g_{S^q} \, .
\end{equation}
Henceforth, denote the metric $\Xi^\ast g_{S^{n+1}}$ by $g$.  Furthermore, the level sets of the coordinate $v$ correspond to generalized Clifford hypersurfaces $C_v$.  The mean curvature of the level set of the coordinate $v$ is given by 
\[
H_v  =  q \,  \cot v  - p  \, \tan v \, .
\]

Next, we introduce canonical coordinates near the point  $\Xi^{-1}(\mu_0) = ( (1, 0, \ldots, 0, ), (1, 0, \ldots, 0), t_\ast)$  in the level set $v = t_\ast$.  (We obtain canonical coordinates in the neighbourhood of the other points of $\Xi^{-1} (\Lambda)$ by symmetry.)  On $S^p \times S^q$, we consider $\bar z \in {\R}^n \longmapsto (\Theta^{(p)}(\bar z), \Theta^{(q)}(\bar z))  \in S^p \times S^q$ to be geodesic normal coordinates near the point $((1, 0,\ldots, 0), (1, 0, \ldots, 0))\in S^p \times S^q$ when $S^p \times S^q$ is endowed with the metric $\cos^2 t_* \, g_{S^p} + \sin^2 t_* \, g_{S^q}$. This metric being a product metric, the geodesic normal coordinates can be defined in such a way that $\bar z = (\bar x , \bar y)$ where $\bar x$ (resp.~$\bar y$) are geodesic normal coordinates close to $(1, 0, \ldots, 0)$ on $S^p$ (resp.~$S^q$) endowed with the metric $\cos^2 t_* \, g_{S^p}$ (resp.~$\sin^2 t_* \, g_{S^q}$).

\subsection{Construction of the Approximate Solution}

To begin the construction of the approximate solution, we first, define the function $\Gamma_\Lambda$ on $C_{t_*}$ that is the unique solution of the equation
\[
{\mathcal L}_{t_*} \, \Gamma_\Lambda = - c_n  \, \sum_{\mu \in \Lambda } \delta_\mu
\]
invariant under the action of $G$. Here, $\delta_\mu$ is the Dirac $\delta$-mass at the point $\mu \in C_{t_*}$ and the constant $c_n$ is the Euclidean volume of $S^{n-1}$.  The first step in our construction is to perturb two generalized Clifford tori on either side of $C_{t_\ast}$ by a proper multiple of the function $\Gamma_\Lambda$ and attach generalized catenoidal necks to the perturbed hypersurface.   This has the effect of reducing the size of error in the mean curvature.  Observe that in geodesic normal coordinates $\bar z$, the function $\Gamma_\Lambda$ can be expanded near $\mu_0$ as 
$$\Gamma_\Lambda  ( \bar z) = 
\begin{cases}
	\tfrac{1}{{n-2}} \, |\bar z |^{2-n} + {\mathcal O}( |\bar z |^{4-n}) &\qquad \mbox{when $n \geq 5$} \\[3mm]
	\tfrac{1}{{2}} \, |\bar z |^{-2} +  {\mathcal O}((\log 1/|\bar z|)) &\qquad \mbox{when $n=4$} \\[3mm]
	|\bar z|^{-1} + \gamma_\Lambda + {\mathcal O}( |\bar z|) &\qquad \mbox{when $n=3$}
\end{cases}$$
where $\gamma_\Lambda \in {\mathbf R}$ is a constant that depends on  $\Lambda$.  By symmetry, this expansion is the same at all other points in $\Lambda$.  For consistency in notation, we agree that $\gamma_\Lambda : = 0$ when $n \geq 4$.

Next, given $t \in (t_*, \frac{\pi}{2})$ we define $t^- \in (0,  t_* )$ {\it via} the relation 
\[
- H_{t^- } = H_t := q \, \cot t  - p \, \tan t \, .
\]
we also define $t^+ :=  t$. Finally, define also the parameter $\eps_t >0$ for $t$ close enough to $t_*$ to be the unique positive solution of 
\[
t^+ - t^- =    \eps_t \,  \int_{-\infty}^\infty \phi^{2-n}(s) \,  ds +  2 \, \eps_t^{n-1} \, \gamma_\Lambda  \, .
\]
where the function $\phi$ has been defined in \eqref{eq:phi}. observe that $t^+ -  t^- =  {\mathcal O} (\eps_t)$.  Finally, set 
\[
r_t := \eps_t^{\frac{n-1}{n}} \, . 
\]
We now define the hypersurface ${\mathcal C}_t^\pm$ to be the image of $C_{t_*} \setminus \bigcup_{\mu \in \Lambda} B_{r_t}(\mu)$ under the mapping 
\[
z \longmapsto \Xi (z , t^+ - \eps_t^{n-1} \, \Gamma_\Lambda (z) )
\]
and also ${\mathcal C}_t^-$ to be the image of $C_{t_*} \setminus \bigcup_{\mu \in \Lambda} B_{r_t}(\mu)$ under the mapping 
\[
z  \longmapsto \Xi (z ,  t^- + \eps_t^{n-1}  \, \Gamma_\Lambda (z) ) \, .
\]
This produces two hypersursurfaces that are close to $C_{t_*} \setminus \bigcup_{\mu \in \Lambda} B_{r_t}(\mu)$ and that have $m_\Lambda$ boundaries. 

We now insert the re-scaled catenoid $\eps_t \, K$ into $S^{n+1}$ by means of the adapted local coordinates as follows. That is, we consider the image of $\{(s, \Theta) \in {\mathbf R}\times S^{n-1} \, : \,  \phi (s) \leq  \eps_t^{-\frac{1}{n}} \}$ under the mapping
\[
(s, \Theta) \longmapsto \Xi \left( z (\eps_t \, \phi (s) \, \Theta),  \tfrac{1}{2}(t^+ +  t^-) + \eps_t \, \psi (s) \right) \, ,
\]
where $\bar z \longmapsto z(\bar z)$ are the geodesic normal coordinates introduced above, along with the images of this hypersurface translated to neighbourhoods of the other points $\mu \in \Lambda$ by the action of the elements of the group $G$. This process produces $m_\Lambda$ hypersurfaces with boundaries, whose union will be denoted by ${\mathcal N}_t$.

The union of the two hypersurfaces ${\mathcal C}_t^\pm$ and ${\mathcal N}_t$ is not a smooth hypersurface; but using cut-off functions we can interpolate between these hypersurfaces in a smooth manner.  This process can be explained as follows.   Because of the invariance under the action of $G$ it is enough to explain how to form the interpolation in the neighbourhood of the point $\mu_0$.  For example when $n \geq 5$,  the graph of $\bar z \longmapsto  t^+ - \eps_t^{n-1} \, \Gamma_\Lambda (\bar z) $ can be expanded near $\mu_0$ in geodesic normal coordinates as 
\[
 t^+ - \eps_t^{n-1} \, \Gamma_\Lambda (\bar z ) =  t^+ -  \frac{_1}{^{n-2}} \, \eps_t^{n-1} \,  |\bar z |^{2-n}  + {\mathcal O} (  \eps_t^{n-1} \, |\bar z |^{4-n} ) \, .
\]
While, changing variables $|\bar z |= \eps_t \, \phi (s) $ with $s >0$, we find with little work that 
\begin{align*}
	\tfrac{1}{2}(t^+ +  t^-) + \eps_t \, \psi(s (\bar z )) &=  \displaystyle \tfrac{1}{2}(t^+ +  t^- ) + \eps_t \, {\displaystyle \int_0^\infty} \phi^{2-n}(v) dv +  \frac{_1}{^{n-2}} \,  \eps_t^{n-1} \, |\bar z |^{2-n}  + {\mathcal O} (\eps_t^{3n-3} \, |\bar z |^{4-3n}) \\
	&=  t^+ +  \tfrac{1}{n-2} \, \eps_t^{n-1} \, |\bar z |^{2-n}  + {\mathcal O} (\eps_t^{3n-3} \, |\bar z |^{4-3n})  \, .
\end{align*}
Observe that, when $|\bar z | \sim r_t$ then both $\eps_t^{n-1} \, |\bar z |^{4-n}$ and $\eps_t^{3n-3} \, |\bar z |^{4-3n}$ are $\mathcal O \big( \eps_t^{\frac{4(n-1)}{n}} \big)$.  This explains why the connected sum is performed when $|\bar z | \sim r_t $ since this precisely minimizes the distance between the graphs of the different summands.

To obtain a smooth hypersurface it is enough to interpolate  between the two graphing functions inside an annulus whose radii are $2 \,  r_t$ and $ r_t /2$.  For example, to interpolate smoothly between the graph of $ t^+ - \eps_t^{n-1} \, \Gamma_\Lambda $ and the graph of $\frac{1}{2} (t^+ +  t^-) + \eps_t \, \psi(s)$ we define the function 
\[
T_t (\bar z ) : =  \eta (\bar z / r_t) \, ( t^+ - \eps_t^{n-1} \, \Gamma_\Lambda (\bar z)) + (1- \eta (\bar z /r_t) )\, \big(\tfrac{1}{2}(t^+ +  t^-) + \eps_t \, \psi (s ( \bar z )) \big)
\]
where $\eta$ is a cut-off function identically equal to $0$ in $B_{1/2}(0)$ and identically equal to $1$ in ${\mathbf R}^n \setminus B_2(0)$.  A similar analysis can be performed for the lower end of the re-scaled catenoid. The final step in the assembly of the different summands of the approximate solution is to extend the above construction so that the resulting surface is invariant under the action of the elements of $G$. We will denote the transition regions by ${\mathcal T}_t$ corresponding to the image of $\bar B_{2 r_t} (0) \setminus B_{r_t/2} (0)$ under the mapping $\bar z \longmapsto \Xi(z(\bar z), T_t(\bar z))$. 

This recipe produces a hypersurface that we will denote $\tilde S_t$, which is a smooth, embedded submanifold of $S^{n+1}$.  It is equal to the connected sum of $\eps_t$-re-scaled catenoids centered at the points of $\Lambda$ and small perturbations of the generalized Clifford tori $C_{t^\pm}$.  Recall that, by construction, these tori have mean curvature equal to $H_t$.  Finally, when $t$ approaches $t_*$, then $\tilde S_t$ approaches two copies of the unique minimal Clifford torus, punctured at the sub-lattice of points $\Lambda$.

The construction of $\approxsol$ in the two lower dimensions $n=3$ and $n=4$ is similar.

\section{The Analysis}
\mylabel{sec:analysis}

\subsection{Deformations of the Approximate Solution}
\mylabel{sec:deform}

The approximate solution $\approxsol$ constructed in the previous section is such that its mean curvature is close to  $H_t$ everywhere except in a small neighbourhood of each gluing point, and it will be shown that it is nevertheless controlled by precise estimates there.  The next task is to set up a means of finding a small deformation of $\approxsol$ whose mean curvature is exactly the constant $H_t$.

To this end, let $\tilde N_t$ be a  choice of unit normal vector field on $\approxsol$ compatible with the orientation.  If $f \in C^{2,\alpha}(\approxsol)$, then  $\exp(f \tilde N_t) (\approxsol)$ is an embedded submanifold of $S^{n+1}$.  The question whether $\exp(f \tilde N_t) (\approxsol)$ has constant mean curvature now becomes a matter of solving a partial differential equation.  We define the \emph{deformation operator} to be  the mapping $\Phi_t : C^{2,\alpha}(\approxsol)  \rightarrow C^{0,\alpha}(\approxsol)$ by 
\[
\Phi_t (f) := H \big( \exp(f \tilde N_t)(\approxsol) \big) ,
\]
where $H(\cdot)$ is the mean curvature operator.  

The deformation operator $\Phi_t$ is a non-linear, partial differential operator on functions $f$ in $C^{2,\alpha}(\approxsol)$ with values in $C^{0,\alpha}$.  The so-called approximate solution $\approxsol$ is an approximation precisely because the estimates of the mean curvature of $\approxsol$ will ensure that $\Phi_t (0) - H_t$ is small (in a suitable norm defined in the next section).  Thus it is hoped that perturbation methods can be used to solve the equation $\Phi_t (f) = H_t$.  The exact formulation of this method is encapsulated in the statement of the \emph{Inverse Function Theorem}.

\begin{nonumthm}[IFT]
    Let $\Phi : \mathcal{B} \rightarrow \mathcal{B}'$ be a smooth map of Banach spaces, set $\Phi(0) = E$ and denote the \emph{linearization} of $\Phi$ at zero by $\mathcal{L} := \Dif \Phi(0)$.  Suppose that $\mathcal{L}$ is bijective and the estimate $\Vert \mathcal{L} X \Vert \geq C \Vert X \Vert$ holds for all $X \in \mathcal{B}$.  Choose $R$ so that if $Y \in \mathcal{B}$ is such that $\Vert y \Vert \leq R$, then $\Vert \mathcal{L} X - \Dif \Phi(Y) X \Vert \leq \tfrac{1}{2} C \Vert X \Vert$.  If $Z \in \mathcal{B}'$ is such that $\Vert Z - E \Vert \leq \frac{1}{2}C R$,  then there exists a unique $X \in \mathcal{B}$ with $\Vert X \Vert \leq R$ so that $\Phi(X) = Z $.  Moreover, $\Vert X \Vert \leq \frac{2}{C} \Vert Z - E \Vert$.
\end{nonumthm}

The first step in applying the IFT to the solution of the problem $\Phi_t (f) = H_t$ is to determine the linearization of $\Phi_t$ at $0$.   We have 
\[
\Dif \Phi_t (0) \, u : = \left. \frac{\dif}{\dif s} \right|_{s=0} \Phi_t ( su)  = \Delta_{\approxsol} u + \Vert \tilde B_t \Vert^2 u + n \, u
\]
where $\Delta_{\approxsol}$ is the Laplacian of $\approxsol$ and $\tilde B_t$ is its second fundamental form.  Henceforth, use the notation $\tilde {\mathcal{L}}_t := \Dif \Phi_t (0)$.

The remaining steps in applying the IFT to the solution of the problem $\Phi_t (f) = H_t$ are the following.   First, appropriate Banach subspaces of $C^{2,\alpha}(\approxsol)$ and $C^{0,\alpha}(\approxsol)$ must be found so that the estimate of $\tilde {\mathcal{L}}_t$ can be achieved.  It must then be shown that $\tilde {\mathcal{L}}_t$ is surjective as a map between these spaces.  Next, estimates in these norms of the non-linear quantities --- the size of $E := \Phi_t (0) - H_t$ and the size of the parameter $R$ giving the variation of $\Dif \Phi_t$ --- must be found.  Note that all these quantities depend {\it a priori} on $t$.    Finally, it must be shown that as a result of these estimates, the quantity $E$ satisfies the inequality $\Vert E \Vert \leq \frac{1}{2} C R$ for all $t$ sufficiently close to $t_*$.  If this holds, then the IFT asserts that a solution of the equation $\Phi_t (f) = H_t$ exists and that it is controlled by the size of $E$ .

\subsection{Function Spaces and Norms}

It does not seem possible to obtain a `good' linear estimate of the form $\Vert \tilde {\mathcal{L}}_t  \, u \Vert \geq C \Vert u \Vert$ with any straightforward choice of Banach subspaces and norms, where `good' in this case means with a constant $C$ independent of $t$.  There are essentially three reasons for this.  The first is that the motion of $\approxsol$ under any isometry of $S^n$ fixes its mean curvature and thus provides an element in the kernel of $\tilde {\mathcal{L}}_t$, also known as a Jacobi field.  Consequently, $\tilde {\mathcal{L}}_t$ is \emph{not} injective on $C^{2,\alpha}(\approxsol)$ due to the Jacobi fields that come from the non-trivial $SO(n)$-rotations of the ambient $S^{n+1}$.   The second reason for the absence of a good linear estimate is that it is possible to perform a motion of $\approxsol$ which consists of an $SO(n+1)$-rotation of only one of the two halves of $\tilde {\mathcal{L}}_t$ while leaving the other half fixed.  The deformation field associated to this motion is equal to the Jacobi field associated to the $SO(n+1)$-rotation on the first half of $\approxsol$, is equal to zero on the other half of $\approxsol$ and interpolates between these two values in the neck regions of $\approxsol$.  This function  approximates  an element of one of the eigenspaces of $\mathcal L_t$ with small eigenvalue.  Thus $\tilde {\mathcal{L}}_t$ possesses small eigenvalues (whose eigenfunctions are called \emph{approximate Jacobi fields}) so that even if one were to choose a Banach subspace of functions transverse to the Jacobi fields coming from isometries of $S^{n+1}$, the constant in the linear estimate would still depend on $t$ in an undesirable manner.  Finally, another source of approximate Jacobi fields is the neck region itself.  It is possible to have a function on $\approxsol$ which is equal to zero away from the neck region and is equal to a Jacobi field of the generalized catenoid within each component of the neck region.  Such an approximate Jacobi field must  `disappear' as $t \rightarrow t_*$ and the necks pinch off, but so long as $t \neq t_*$, these functions contribute to the size of the constant $C$ in the linear estimate for $\tilde {\mathcal{L}}_t$.

The three problems listed above will be dealt with here in two ways. First, the symmetry group $G$ of the approximate solution must be exploited.  It turns out that the Jacobi fields, both approximate and true, do \emph{not} share these same symmetries; thus working in a space of functions possessing these symmetries will rule out the existence of small eigenvalues. 
Indeed, under the assumption of Theorem~\ref{th:1} and according to the result of Corollary~\ref{cor:nojacobi} there is no non-trivial solution of ${\mathcal L}_{t_*} \, u = 0$ that is invariant under the action of $G$.  This means that restricted to the set of $G$-invariant functions, the operator ${\mathcal L}_{t_*} :  { C}^{2,\alpha} (C_{t_*}) \longrightarrow { C}^{0, \alpha}(C_{t_*})$ is an isomorphism.  Second, it is necessary to use a somewhat non-standard norm to measure the `size' of functions $f \in C^{2, \alpha}(\approxsol)$ in order to properly determine the dependence on the parameter $t$ of the various estimates needed for the application of the Inverse Function Theorem.  A \emph{weighted} Schauder norm will be used for this purpose, defined \emph{via} a \emph{weight function}.  As usual, one can say without loss of generality that the weight function is invariant with respect to the symmetry group $G$.

First let $r<1$ be some fixed radius that is determined by the following two requirements: $r$ is such that adapted local coordinates can be defined inside $B_r(\mu)$ for each $\mu \in \Lambda$; and $B_{2r}(\mu)$ and $B_{2r}(\mu')$ are disjoint for all $\mu \neq \mu' \in \Lambda$.   

\begin{defn}
    \mylabel{defn:weight}
    The \emph{weight function} $\zeta_t : \approxsol \rightarrow \R$ is defined by
\begin{equation*}
    \zeta_t (z, v) =
    \begin{cases}
    	1 &\qquad \Xi (z,v) \in \mathcal C_t^\pm \setminus \bigcup_{\mu \in \Lambda} B_{2r}(\mu) \\
        \mbox{Interpolation} &\qquad \Xi (z,v) \in \mathcal C_t^\pm \cap \big[ \bar B_{2r}(\mu) \setminus B_{r}(\mu) \big] \: \mbox{for some} \: \mu \in \Lambda \\
        |\bar z |  &\qquad  \Xi (z(\bar z), v) \in {\mathcal C}_t^\pm \cap B_{r} (\mu)  \: \mbox{for some} \: \mu \in \Lambda \\
                \mbox{Interpolation} & \qquad \Xi (z, v) \in {\mathcal T}_t \\
        \eps_t \, \cosh s & \qquad  \Xi (z( \eps_t \phi(s) \Theta ), \frac{1}{2} (t^++t^-)+ \eps_t \psi(s) ) \in {\mathcal N}_t \, .
    \end{cases}
\end{equation*}
\end{defn}

\noindent Next, let $T$ be any tensor on $\approxsol$, and recall the notation
\begin{equation*}
    \Vert T \Vert_{0,\approxsol} = \sup_{p \in \approxsol} \Vert T(p) \Vert \\
    \qquad \mbox{and} \qquad [T]_{\alpha,\approxsol} = \sup_{p, p' \in \approxsol} \frac{\Vert T(p) - \Pi (T(p'))  \Vert}{\mathrm{dist}(p, p')^\alpha} \, ,
\end{equation*}
where the norms and the distance function that appear are taken with respect to the induced metric of $\approxsol$, while $\Pi$ is the parallel transport operator from $p$ to $p'$ with respect to this metric.  Now define
\begin{equation}
    \mylabel{eqn:weightnorm}
    \vert f \vert_{\ckbg(\approxsol)} := \vert \zeta_t^{-\gamma} f \vert_{0,\approxsol} + \Vert \zeta_t^{-\gamma + 1} \nabla f \Vert_{0,\approxsol} + \cdots + \Vert \zeta_t^{-\gamma + k} \nabla^k f \Vert_{0,\approxsol} + [ \zeta_t^{-\gamma + k + \alpha} \nabla^k f ]_{\alpha, \approxsol} \, .
\end{equation}
Again, the norms and derivatives which appear here are taken with respect to the induced metric of $\approxsol$. It is easy to check that the space of $C^{k,\alpha}$ functions on $\approxsol$ measured with respect to the norm \eqref{eqn:weightnorm} is a Banach space.  Henceforth, denote this space of functions by $\ckbg(\approxsol)$.

A solution of the deformation problem will be found in a space of $C^{2,\alpha}$ functions on $\approxsol$.  To do so, it will be necessary to insist that the functions in this space inherit the symmetries of $\approxsol$ since this will have the effect of ruling out the existence of  the Jacobi fields and the approximate Jacobi fields of $\approxsol$ which are the analytic obstructions preventing the inversion of the deformation operator.  The following space will meet these needs.

\begin{defn}
    Let $\mathcal{B}^{k,\alpha}_\gamma (\approxsol) : =  \{ f \in C^{k,\alpha}_\gamma(\approxsol) \: : \: f \circ \sigma = f \: \mbox{for all} \: \sigma \in G \}$.
\end{defn}

 Clearly, the operator $\Phi_t$ is a well-defined map from $\Btbg$ to $\Bobg$ that is smooth in the Banach space sense.  The linearized operator $\tilde{\mathcal{L}}_t : \Btbg \rightarrow \Bobg$ is bounded and satisfies
    \begin{equation*}
        \vert \tilde{\mathcal{L}}_t \, u  \vert_{\cobg(\approxsol)} \leq C \,\vert u \vert_{\ctbg(\approxsol)} 
    \end{equation*}
    where $C$ is independent of $t$, chosen close enough to $t_*$.  Finally, for any $\gamma \in \R$, there exists another constant $C$ independent of $t$ so that $\tilde{\mathcal{L}}_t$ satisfies the elliptic estimate
    \begin{equation}
        \mylabel{eqn:ellest}
        \vert u \vert_{\ctbg(\approxsol)} \leq C \big( \vert \tilde{\mathcal{L}}_t \, u \vert_{\cobg(\approxsol)} + \vert \zeta_t^{-\gamma} u \vert_{0, \approxsol} \big) \, .
    \end{equation}
This follows at once from Schauder's estimates applied on the different summands constituting $\tilde S_t$. 

\subsection{The Linear Estimate}

\newcommand{\ctbgi}{C^{2,\alpha}_\gamma(\tilde S_{t_i})}
\newcommand{\cobgi}{C^{0,\alpha}_{\gamma - 2}(\tilde S_{t_i})}
\newcommand{\approxsoli}{\tilde S_{t_i}}
\newcommand{\Btbgi}{\mathcal{B}^{2,\alpha}_{\gamma} (\tilde S_{t_i})}

The most important estimate needed to solve the equation $\Phi_t (f) = H_t$ by means of the Inverse Function Theorem is the estimate from below of the linearization $\tilde{\mathcal{L}}_t$ by a constant independent of $t$.  The purpose of this section is to prove this estimate using an argument by  contradiction, in which it is assumed that such a lower bound does not exist.

\begin{prop}
    \mylabel{thm:linest}
Suppose $2-n < \gamma < 0$.  Then the linearized operator $\tilde{\mathcal{L}}_t : \Btbg \rightarrow \Bobg$ satisfies
\[
\vert \tilde{\mathcal{L}}_t \, u  \vert_{\cobg(\approxsol)} \geq C \vert u \vert_{\ctbg(\approxsol)}
\]
where $C$ is a constant independent of $t$ close enough to $t_*$.
\end{prop}

\begin{proof}
Observe that Schauder's elliptic estimates  imply that it is enough to prove that 
\[
\vert \zeta^{2-\gamma}_t \, \tilde{\mathcal{L}}_t \, u  \vert_{L^\infty (\approxsol)} \geq C \vert \zeta^{-\gamma}_t \,  u \vert_{L^\infty (\approxsol)}
\]
where $C$ is a constant independent of $t$ close enough to $t_*$.

We argue by contradiction.  Suppose that $t_i \rightarrow t_*$ and that there is a sequence of functions $u_i$ defined on $\approxsoli$ that are invariant under the action of $G$ along with a sequence of linearized operators $\tilde{\mathcal{L}}_{t_i}$ satisfying the following estimates:
\[
\lim_{i\rightarrow \infty} \vert \zeta^{2-\gamma}_{t_i} \, \tilde {\mathcal{L}}_{t_i} u_i \vert_{L^\infty (\approxsoli )} = 0  \qquad \mbox{and} \qquad \vert \zeta^{-\gamma}_{t_i} \, u_i \vert_{L^\infty (\approxsoli )} = 1 \, .
\]
Moreover, one can assume that $\approxsoli$ converges in a smooth enough sense to two copies of the unique minimal generalized Clifford hypersurface with the gluing points removed (denote this hypersurface by $C_{t_*}\setminus \Lambda$) and that the operators $\tilde{\mathcal{L}}_{t_i}$ converge to the Jacobi operator on $C_{t_*}$, which is simply $\mathcal{L}_{{t_*}}$.  Let $q_i$ be a point where $\big( \zeta_{t_i}(q_i) \big)^{-\gamma} \, | u_i (q_i) | = 1$; then up to a subsequence, either $q_i \rightarrow q \in C_{t_*}\setminus \Lambda$, or else $q_i$ converges to a point of $\Lambda$.  These two scenarios will be ruled out in turn.  In what follows, adapted local coordinates will always be used in the neighbourhood of $q_i$.

\medskip \noindent \scshape Case 1. \upshape  Suppose $q_i \rightarrow q \in C_{t_*}\setminus \Lambda$.

\medskip

In this case, one uses elliptic estimates together with Arzela-Ascoli's theorem, to prove that (up to a subsequence) $u_i \rightarrow u_\infty$ uniformly on compact subsets of $C_{t_*}\setminus \Lambda$.  The limit function $u_\infty$ satisfies $\mathcal{L}_{{t_*}} \, u_\infty = 0$ on $C_{t_*}\setminus \Lambda$ and 
\begin{equation}
| (\mbox{dist} (\cdot , \Lambda) )^{-\gamma} \, u_\infty |_{L^\infty (C_{t_*})} = 1.
\label{eq:fdfd}
\end{equation}
Finally, $u_\infty$ is invariant under the action of $G$.  Since we have assumed that $\gamma > 2-n$, the singularities are removable and hence $u_\infty$ is smooth. But by assumption, no nontrivial element of the kernel of ${\mathcal L}_{t_*}$ is invariant under the action of $G$.   This is clearly in contradiction with (\ref{eq:fdfd}) and hence rules out Case 1.

\medskip \noindent \scshape Case 2. \upshape  Suppose without loss of generality that $q_i$ converges to the gluing point $\mu_0 \in \Lambda$.

\medskip

This second case divides into two subcases. First assume that (up to a subsequence)  $q_i$ belongs to ${\mathcal N}_{t_i}$ so that it can be written as $$q_i = \Xi \big( z ( \eps_{t_i} \, \phi(s_i) \Theta_i ), \tfrac{1}{2}( t_i^+ + t_i^-) +   \eps_{t_i} \, \psi (s_i) \big) \, .$$ Further assume that the (up to a subsequence)  the sequence $s_i$ is bounded and even converges to $s_\infty$.  The use of elliptic estimates together with Arzela-Ascoli's theorem is enough to prove that (up to a subsequence) $u_i \rightarrow u_\infty$ uniformly on compact subsets of $K$.  The limit function $u_\infty$ satisfies $\mathcal{L}_{K} \, u_\infty = 0$ on $K$ and 
\begin{equation}
| (\cosh s)^{-\gamma} \, u_\infty |_{L^\infty (K)} = 1.
\label{eq:fdfd2}
\end{equation}
Finally, since $u_i$ is invariant under the action of $\rho \in G$, the limit function $u_\infty$ is invariant under the action of the symmetry with respect to the $x_{n+1}$-axis.  Since we have assumed that $\gamma <0$, then the result of Proposition~\ref{prop:catjacobi} implies that $u_\infty =0$ which is clearly in contradiction with (\ref{eq:fdfd2}). This rules out this first subcase.

Now it remains to consider the case that is not covered by the first subcase. This time  $q_j$ converges to $\mu_0$ but at a slower rate and  the use of elliptic estimates together with Arzela-Ascoli's theorem is enough to prove that (up to a subsequence) $u_i \rightarrow u_\infty$ uniformly on compact subsets of ${\mathbb R}^p \times {\mathbb R}^q \setminus \{0,0\}$, where this space is endowed with the metric $g_* : = (\cos t_*)^2 \, \mathring g_p + (\sin t_*)^2 \, \mathring g_q$ where $\mathring g_N$ denotes the Euclidean metric on $\R^N$.  The limit function $u_\infty$ satisfies $\Delta_{g_*} \, u_\infty = 0$ on ${\mathbb R}^p \times {\mathbb R}^q \setminus \{0,0\}$ and 
\begin{equation}
\big| |\bar z |^{-\gamma} \, u_\infty \big|_{L^\infty } = 1.
\label{eq:fdfd3}
\end{equation}
Since we have assumed that $2-n <\gamma <0$, then this clearly implies that $u_\infty =0$ which is clearly in contradiction with (\ref{eq:fdfd3}). This rules out this second and last subcase.

Having ruled out all possible cases, the proof of the claim and hence the proof of the result are complete. 
\end{proof}

\subsection{The Estimate of the Mean Curvature of the Approximate Solution}

As mentioned earlier, the proof of Theorem~\ref{th:1} requires two more estimates in addition to the one from the previous section.  The first of these it so show that $\Phi_t (0) - H_t$ is small in the $\cobg$ norm. The following calculations are generalizations of those carried out in \cite{mepacard1}.

\begin{prop}
    \mylabel{prop:error}
    The quantity $\Phi_t (0)$, which is the mean curvature of $\approxsol$, satisfies the following estimate. Assume that $\gamma  > 2-n $ is fixed. If $t$ is sufficiently close to $t_*$, then there exists a constant $C$ independent of $t$ so that
    \begin{equation}
        \mylabel{eqn:error}
        \big\vert \Phi_t (0) - H_t \big\vert_{\cobg (\tilde S_t)} \leq C \, \eps_t^{2-\gamma} \, .
    \end{equation}
   \end{prop}
\begin{proof}
To estimate the mean curvature in the regions ${\mathcal  T}_t$ and $\mathcal C_t^\pm$, we first compute the mean curvature of the graph of a function $u : S^p \times S^q \rightarrow \R$ parameterized by $z  \longmapsto \Xi  ( z , \, u (z))$ where $\Xi$ is the toroidal coordinate embedding into $S^{n+1}$.  Henceforth we will use the notation that a comma denotes partial differentiation, such as $u_{,i} = \partial_{z_i} u$, and repeated indices are summed. The tangent vectors of this surface are given by $T_j : =  \partial_{z_j } + u_{,j} \, \partial_v$ and it is easy to check that the induced metric is given by
\[
\bar g =  du \otimes du  + \cos^2 u \, g_{S^p} + \sin^2 u \, g_{S^q}.
\]
The normal vector field $N$ can be written as $N : = \bar N / ||Ê\bar N||$ where  $\bar N : = \partial_t - a^j \, T_j$ and the coefficients $a^j$ are determined so that $\bar N$ is normal to the surface. One finds the explicit expressions $a^j =  \bar g^{jk} u_{,k}$ and $||Ê\bar N||^2 =  1+ \bar g^{jk} u_{,j} u_{,k}$.   We now compute
\begin{align*}
2 \, g (\nabla_{T_i} T_j,  \bar N) & = - g (\nabla_{T_i} \bar N , T_j)  - g(\nabla_{T_j} \bar N , T_i) \\
& = - ( g(\nabla_{T_j} \partial_v , T_i)  + g(\nabla_{T_i} \partial_v , T_j) ) + (\partial_{z_i} g(a^k T_k, T_j) + \partial_{z_j} g(a^k T_k, T_i))  -2 \, a^k \, g(T_k, \nabla_{T_i} T_j)\\
& = - ( g(\nabla_{T_j} \partial_v, T_i)  + g(\nabla_{T_i} \partial_v, T_j) ) + (\partial_{z_i} g(\partial_v , T_j) + \partial_{z_j} g(\partial_v , T_i))  - 2 \, a^k \, g(T_k, \nabla_{T_i} T_j)\\
& = - ( g(\nabla_{T_j} \partial_v , T_i)  + g(\nabla_{T_i} \partial_v , T_j) ) + 2 \, u_{,ij}  - 2 \bar \Gamma_{ij}^k \, u_{,k}
\end{align*}
where $\bar \Gamma_{ij}^k = \frac{1}{2} \bar g^{l k} \, ( \bar g_{jl,i}+ \bar g_{il,j} - \bar g_{ij,l})$ are the Christoffel symbols of $\bar g$.  To evaluate the first terms, we consider the parameterization of a neighbourhood of $C_{t_*}$ given by $\tilde \Xi  : (z, \xi) \longmapsto \Xi (z, u(z)+\xi)$ so that $\partial_v = \tilde \Xi_\ast \partial_\xi$. Therefore, we can write 
\[
 g(\nabla_{T_j} \partial_v , T_i)  + g(\nabla_{T_i} \partial_v , T_j)  =  \partial_\xi \, g( T_i, T_j) \big|_{\xi=0} =  2 \, \cos u \, \sin u \,  \big[ g_{S^q} - g_{S^p}\big]_{ij} \, .
\]
Collecting these, we  obtain the second fundamental form
\begin{equation*} 
 \| \bar N\| \, \bar B =  \cos u \,  \sin u  \, ( g_{S^q} - g_{S^p} ) +  \big( u_{,ij} - \bar \Gamma^k_{ij} u_{,k} \big)  d\bar z_i\, d\bar z_j \, .
\end{equation*}
Finally, we get the mean curvature by taking the trace of $\bar B$ with respect to $\bar g$.

We now specialize this computation to the case where  
\[
u=u_t :=  t^+ - \eps_t^{n-1} \, \Gamma_\Lambda \qquad \mbox{and}Ê\qquad \mathrm{dist}(z, \Lambda ) \geq r_t
\] 
so that we obtain the mean curvature of ${\mathcal C}^+_t$.  We estimate the metric coefficients Christoffel symbols as
\begin{align*}
	\bar g &= \cos^2 u_t  \, g_{S^p}  + \sin^2 u_t \, g_{S^q}  + \mathcal O (\eps_t^{2n-2} |\bar z|^{2-2n}) \\[1ex]
	\bar \Gamma_{ij}^k &= 
	\begin{cases}
		\frac{1}{\cos^2 u_t} \stackrel{\mbox{\tiny$p$}\: \:\:}{\Gamma_{ij}^k}  + \mathcal O (\eps_t^{n-1} |\bar z|^{1-n})  &\quad \mbox{if $i,j$ and $k$ refer to $S^p$ terms}\\
		\frac{1}{\sin^2 u_t} \stackrel{\mbox{\tiny$q$}\: \:\:}{\Gamma_{ij}^k} + \mathcal O (\eps_t^{n-1} |\bar z|^{1-n})  &\quad \mbox{if $i,j$ and $k$ refer to $S^q$ terms}\\
		\mathcal O(\eps_t^{n-1} |\bar z |^{1-n})  &\quad \mbox{if $i,j$ and $k$ mix $S^p$ with $S^q$ terms}
	\end{cases}
\end{align*}
where  $\stackrel{\mbox{\tiny$p$}\: \:\:}{\Gamma_{ij}^k} $ (resp.~$\stackrel{\mbox{\tiny$q$} \:\:\:}{\Gamma_{ij}^k} $) are the Christoffel symbols of $g_{S^p}$ (resp.~$g_{S^q}$).  For example, when $i,j$ and $k$ refer to $S^p$ terms, we have 
\[
\bar \Gamma_{ij}^k =\tfrac{1}{\cos^2 u_t} \stackrel{\mbox{\tiny$p$}\: \:\:}{\Gamma_{ij}^k}  + \mathcal O (\eps_t^{n-1} |\bar z|^{1-n}) + \mathcal O (\eps_t^{2n-2} |\bar z|^{1-2n} ) \, .
\]
But since we are working in the range where $ \mathrm{dist}(z, \Lambda ) \geq r_t$, we can use the fact that
\[
\eps_t^{2n-2} |\bar z|^{1-2n}  = \mathcal O (\eps_t^{n-1} |\bar z|^{1-n}) 
\]
so that the estimate simplifies into 
\[
\bar \Gamma_{ij}^k =\frac{1}{\cos^2 u_t} \stackrel{\mbox{\tiny$p$}\: \:\:}{\Gamma_{ij}^k}  + \mathcal O (\eps_t^{n-1} |\bar z|^{1-n}) \, .
\]

We conclude that for $ \mathrm{dist}(z, \Lambda ) \geq r_t$, 
\[
\| \bar N\| \, \bar B =    \cos u \,  \sin u  \, ( g_{S^q} - g_{S^p} ) +  \big( \nabla_{ij}^\ast u + \mathcal O (\eps_t^{2n-2} \, |\bar z|^{2-2n}) \big)  d\bar z_i\, d\bar z_j  \, ,
\]
where $\nabla^\ast_{ij}$ is the covariant derivative of $S^p$ (resp.~$S^q$) if $i,j$ refer to $S^p$ terms (resp.~$S^q$ terms).  Consequently, 
\begin{align*}
\| \bar N\| \, H &= \tfrac{1}{\cos^2 u_t} \, \Delta_{S^p} u_t+ \tfrac{1}{\sin^2 u_t} \, \Delta_{S^q} u_t  +  (q \, \cot u_t  - p \, \tan u_t ) + \mathcal O(\eps_t^{2n-2} |\bar z |^{2-2n} ) \\[1ex]
	&=  (q \, \cot t^+  - p \, \tan t^+ ) - \left(  \Delta_{C_{t^+}} + \left( \tfrac{p}{\cos^2 t^+} + \tfrac{q}{\sin^2 t^+}\right) \right)  \eps_t^{n-1} \Gamma_\Lambda  + \mathcal O(\eps_t^{2n-2} |\bar z |^{2-2n}) + \mathcal O ( \\[1ex]
	&= H_t -  \mathcal L_{t} ( \eps_t^{n-1} \Gamma_\Lambda ) +  \mathcal O(\eps_t^{2n-2} |\bar z |^{2-2n})  \\
	&= H_t + \mathcal O (\eps_t^n \, |\bar z |^{-n}) +  \mathcal O(\eps_t^{2n-2} |\bar z |^{2-2n})
\end{align*}
where $\Delta_{C_t}$ is the Laplace operator of $C_t$ in which case we have used the fact  that $\mathcal L_{t_*} \, \Gamma_\Lambda = 0$ away from the point $\mu_0$ as well as the fact that $\mathcal L_{t_*}- \mathcal L_{t}$ is a second order differential operator whose coefficients are  bounded by a constant times $\eps_t$.  We have also used the various fall-off behaviors of $\Gamma_\Lambda$ and its derivatives to obtain the result above.  Observe that $\eps_t^{2n-2} |\bar z|^{2-2n}  = \mathcal O (\eps_t^{n} |\bar z|^{-n}) $ so that we finally get 
\[
H = H_t + \mathcal O (\eps_t^{n} |\bar z|^{-n}) 
\]
by taking the trace with respect to $\bar g$.

The corresponding estimates in ${\mathcal C}^-_t$ and also in ${\mathcal T}_t$ are obtained using similar computations. Observe that the cut-off functions used in $\mathcal TÑ_t$ induce another discrepancy that can be estimated by a constant times $\eps_t^{n-1}Ê\, r_t^{2-n} = \mathcal O (\eps_t^{n} r_t^{-n}) $ when $n \neq 4$ (and by a constant times  $\eps_t^{3}Ê\, r_t^{-2} \, (\log 1/r_t ) = \mathcal O (\eps_t^{4} r_t^{-4})$ in dimension $n=4$). In any case $|\bar z| \sim r_t$ in $\mathcal T_t$ and hence we still have $H = H_t +  \mathcal O (\eps_t^{n} |\bar z|^{-n})$ in the transition region $\mathcal T_t$. 

It remains to compute the mean curvature of the neck region ${\mathcal N}_t$. Since the center of the neck region is not a graph over the level sets of constant $t$, the previous calculation does not help us.  Thus we compute directly the mean curvature of the surface parameterized by  
\begin{equation}
\label{eqn:catemb}
(s, \nu) \longmapsto \Xi \big( z(\eps_t \phi(s)  \Theta (\nu)) , \tfrac{1}{2}(t^+ +  t^-) + \eps_t  \, \psi (s) \big)
\end{equation}
where $s \in \R$ satisfies $\phi(s) \leq \eps_t^{-\frac{1}{n}}$ and $\nu \in S^{n-1} \longmapsto  \Theta(\nu) \in \R^n$ is a parametrization of the unit $(n-1)$-sphere.  Since the embedding \eqref{eqn:catemb} is defined using geodesic normal coordinates $(\bar z, v)$, the background metric is of the form
\begin{align*}
	g &= dv^2 + \frac{\cos^2 v}{\cos^2 t_\ast} \sum_{i,j = 1}^p \big( \delta_{ij} + Q^p_{ij} \big) d \bar x_i \, d \bar x_j + \frac{\sin^2 v}{\sin^2 t_\ast} \sum_{i,j = p+1}^n \big( \delta_{ij} + Q^q_{ij} \big) d \bar y_i \, d \bar y_j 
\end{align*}
where the components $Q^p_{ij}$ and $Q^q_{ij}$ satisfy 
\[
|Q| + |\bar z| \, |DQ| + |\bar z|^2 |D^2Q| = \mathcal O(|\bar z|^2).
\]

Now we compute the tangent vectors, the induced metric and the normal vector of the region $\mathcal N_t$ in these coordinates.  Use Greek letters to refer to the components of $\nu$, such as $\Theta^i_{,\alpha} = \frac{\partial \Theta^i}{\partial \nu_\alpha}$, and use a dot to indicate differentiation with respect to $s$.  The tangent vectors are given by 
\[
T_\alpha := \eps_t \, \phi \,  \Theta^i_{,\alpha} \, \partial_{\bar z_i}  \qquad \mbox{and} \qquad T_s : =  \eps_t \, \partial_s \phi  \,  \Theta^i \, \partial_{\bar z_i} + \eps_t \, \phi^{2-n} \, \partial_v .
\]
The induced metric is then given by
\begin{align*}
	\bar g &=  \eps_t^2 \phi^2 \big( d s^2 + g_{S^{n-1}} \big) + Q ds^2 + \sum_{\alpha = 1}^{n-1} Q_{\alpha } d\nu_\alpha \, ds + \sum_{\alpha,\beta = 1}^{n-1} Q_{\alpha \beta} d\nu_\alpha \, d\nu_\beta 
\end{align*}
where this time the components $Q, Q_{\alpha}, Q_{\alpha \beta}$ satisfy the estimate 
\[
|Q| + |DQ| + |D^2Q| = \mathcal O(\eps_t^3 \cosh^2 s)  + \mathcal O ( \eps_t^4 \, \cosh^4 s) .
\] 
The normal vector field $N$ can be written as $N : = \bar N / ||Ê\bar N||$ with $\bar N : = N_0 - a^s \, T_s - a^\alpha \, T_\alpha$ where 
$$
N_0 : =  - \phi^{1-n} \,  \Theta^i \, \partial_{\bar z_i}  + \partial_s \log \phi  \, \partial_vÊ
$$
and where the coefficients $a^s$, $a^ \alpha $ are determined so that $\bar N$ is normal to the surface.  Using the fact that we are only interested in the region where $\phi(s) \leq \eps_t^{-\frac{1}{n}}$, one finds the estimates
\begin{align*}
	a^s &=  {\mathcal O}(  \cosh^{-n} s) + \mathcal O ( \eps_t \, \cosh^{2-n} s)\\
	a^\alpha &=  {\mathcal O}(  \cosh^{-n} s) + \mathcal O ( \eps_t \, \cosh^{2-n} s) \\
	||Ê\bar N||^2 &=  1 +  {\mathcal O}(\eps_t \, \cosh^{2-2n} s) + \mathcal O ( \eps_t^2 \, \cosh^{4-2n} s) \, .
\end{align*}

We compute as above
\begin{align*}
	 2 \, \|\bar N \| \, B(T_i, T_j) = & - g(\nabla_{T_i} \bar N, T_j) - g(\nabla_{T_j} \bar N, T_i) \\
	 = & - (g(\nabla_{T_i} N_0, T_j) + g(\nabla_{T_j} N_0, T_i) ) + (g(N_0, T_i)_{,j} + g(N_0, T_j)_{,i} ) \\
	 & -  a^k \, ( \bar g_{jk,i} + \bar g_{ik,i} -  \bar g_{ij,k}) 
\end{align*}
where $i, j, k$ can be $s$ or $\alpha$.  It is easy to check that \[
g(N_0, T_i)_{,j} + g(N_0, T_j)_{,i}  = \mathcal O(\eps_t^2 \, \cosh^{2-n} s) + \mathcal O(\eps_t^3 \, \cosh^{4-n} s)
\] 
and also that 
\[
a^k \, ( \bar g_{jk,i} + \bar g_{ik,i} -  \bar g_{ij,k})  = \mathcal O(\eps_t^2 \, \cosh^{2-n} s) +  \mathcal O(\eps_t^3 \, \cosh^{4-n} s).
\]  

We now compute the first terms.  As above, our calculations are simplified by considering the local parameterization 
\[
\hat \Xi : (s, \nu, \xi) \longmapsto \Xi \left( z \big( (\eps_t \, \phi( s)  - \xi \,  \phi^{1-n} (s) ) \, \Theta(\nu) \big),  \tfrac{1}{2} (t^+ + t^-) + \eps_t  \psi(s)  +  \xi \, \partial_s \log \phi (s) \right)
\]
so that $N_0  = \hat \Xi_\ast \, \partial_\xi$ and hence $g(\nabla_{T_i} N_0, T_j) + g(\nabla_{T_j} N_0, T_i) = \partial_\xi  g(T_i,T_j) \big|_{ \xi=0}$.  We obtain with little work
\[
\partial_\xi g(T_i, T_j) \big|_{\xi=0} = 
\begin{cases}
	 	 2 \eps_t  \, (n-1) \phi^{2-n} + 2\,  \eps_t^2   \, (\partial_s \phi )^2 \big(  \cot t_* d\bar y^2 - \tan t_* d\bar x^2 \big) (\Theta, \Theta) \\+ {\mathcal O}(\eps_t^2 \, \cosh^{2-n} s)  + {\mathcal O}(\eps_t^3 \, \cosh^{4-n} s) + \mathcal O (\eps_t^3 \, \cosh^2 s) & \quad \mbox{when $i=j=s$} \\[3mm]
	 2\,  \eps_t^2   \, \partial_s \phi \, \phi \, \big(  \cot t_* d\bar y^2 - \tan t_* d\bar x^2 \big) (\Theta, \Theta_{,\alpha})  \\+ {\mathcal O}(\eps_t^2 \, \cosh^{2-n} s) 
	+ {\mathcal O}(\eps_t^3 \, \cosh^{4-n} s)+ \mathcal O (\eps_t^3 \, \cosh^2 s)& \quad \mbox{when $i=s, j= \alpha$} \\[3mm]
	 	 - 2 \eps_t  \, \phi^{2-n} \big[g_{S^{n-1}}]_{\alpha \beta}  + 2\,  \eps_t^2   \,\phi^2 \big(  \cot t_* d\bar y^2 - \tan t_* d\bar x^2  \big) (\Theta_{, \alpha}, \Theta_{,\beta})  \\+ {\mathcal O}(\eps_t^2 \, \cosh^{2-n} s) + {\mathcal O}(\eps_t^3 \, \cosh^{4-n} s) + \mathcal O (\eps_t^3 \, \cosh^2 s) & \quad \mbox{when $i=\alpha, j= \beta$}
\end{cases}
\]
When we take the trace of the second fundamental form computed above, the leading-order terms coming from $\partial_\xi g(T_i, T_j) \big|_{\xi=0}$ vanish because $(s, \nu) \longmapsto (\phi(s) \Theta(\nu), \psi(s))$ is a minimal embedding into $\R^{n+1}$. The next-leading-order terms also vanish because  $C_{t_*}$ is minimal in $S^{n+1}$.  Thus only the remaining terms contribute to the estimate of the mean curvature.  Since $H_t = \mathcal O (\eps_t)$, we conclude that
\[
H - H_t = \mathcal O (\eps_t) + {\mathcal O}(\cosh^{-n} s) + \mathcal O (\eps_t \, \cosh^{2-n} s) =  {\mathcal O}(\cosh^{-n} s) \, .
\]
using the fact that $\phi \leq \eps_t^{- \frac{1}{n}}$.

It remains to collect the estimates in all the various regions of $\approxsol$ and perform the estimate in the weighted H\"older norm.  We get 
\begin{align*} 
	| H(\tilde S_t) - H_t|_{C^0_{\gamma-2} (\approxsol)} &\leq \sup_{z \in \approxsol} \big| \zeta^{2-\gamma}(z) \cdot \big( H(z) - H_t \big) \big| \\
	&\leq | C \, ( \eps_t^n + \eps_t^{2-\gamma}) \\
	&\leq C \, \eps_t^{2-\gamma}
\end{align*}
for some constant $C$ independent of $t$, close enough to $t_*$. Here we have used the fact that $\gamma >2-n$. This completes the estimate of the mean curvature. The estimate of its H\"older coefficient follows similarly. 
\end{proof}

\subsection{The Nonlinear Estimate}

The remaining estimate that is needed to invoke the Inverse Function Theorem is  to show that $\Dif \Phi_t (f)  - \tilde{\mathcal{L}}_t $ can be made to have small operator norm if $f$ is chosen sufficiently small in the $\ctbg  (\approxsol )$ norm.  Once these estimates are given, it will be possible to conclude the proof of Theorem~\ref{th:1} by invoking the Inverse Function Theorem.  

\begin{prop}
    \mylabel{prop:nonlin}
Given $\kappa >0$, there exists $C_\kappa >0$ such that, for all $t$ close enough to $t_*$ and for all $\vert f \vert_{\ctbg  (\approxsol )} \leq \kappa \, \eps_t^{2-\gamma}$ then
    \begin{equation}
        \mylabel{eqn:nonlin}
        \big\vert \Dif \Phi_t (f) \, u - \tilde{\mathcal{L}}_t \, u \big\vert_{\cobg (\approxsol )} \leq C_\kappa \, \eps_tÊ\,  \vert u \vert_{\ctbg  (\approxsol )}
    \end{equation}
     for any $u \in \ctbg(\approxsol)$.
\end{prop}

\begin{proof} To begin with, we consider a hypersurface $\Sigma$ embedded in a Riemannian manifold and $\Sigma_f$ normal perturbation of it for some small function $f$. We assume that $|f|_{{\mathcal C}^{2, \alpha}} \leq c$ where $c$ is some small constant. It is clear that the difference between the Jacobi operator about $\Sigma$ and $\Sigma_f$ is a second order differential operator whose coefficients are bounded in ${C}^{0, \alpha}$ topology by a constant times  $| f |_{{C}^{2, \alpha}}$. 

Now we consider a point $p \in \approxsol$ and a geodesic ball centered at $p$ of radius $r \sim  \zeta_t (p)$. We consider the normal graph over this geodesic ball for a function $f$ whose ${C}^{2, \alpha}_\gamma (\approxsol)$ norm is controlled by  $\kappa \, \eps_t^{2-\gamma} $.  If we blow up $S^{n+1}$ by a factor $1/\zeta_t (p)$. We now have a geodesic ball $D$ of radius $\sim 1$ on the dilated hypersurface in $S^{n+1} (1/ \zeta_t (p))$ and a normal graph over this ball for a function whose ${C}^{2, \alpha}$ norm is controlled by  $\kappa \, \eps_t^{2-\gamma} \, (\zeta_t (p))^{\gamma-1}$. The ball $D$ depends on $t$ but its geometry is controlled uniformly as $t$ tends to $0$, which is a consequence of the definition of the weight function. We can apply the above argument to check that the difference between the Jacobi operator about $D$ and its normal perturbation is a second order differential operator whose coefficients are bounded in ${C}^{0, \alpha}$ topology by a constant times $\kappa \, \eps_t^{2- \gamma} \, (\zeta_t (p))^{\gamma-1}$. 

By performing the dilation backward, we obtain that the difference of the Jacobi operators between $\approxsol$ and its normal perturbation. Observe that the backward dilation multiplies the result by $1/(\zeta_t(p))^2$ but this coefficient is absorbed by the fact that the norm on the left hand side of (\ref{eqn:nonlin}) involves $\gamma-2$ and not $\gamma$.  Varying $p$ along $\approxsol$ we find that the worst estimate occurs precisely in the neck when $\zeta_t(p) \sim \eps_t$. This implies readily that $\big\vert \Dif \Phi_t (f) \, u - \tilde{\mathcal{L}}_t \, u \big\vert_{\cobg (\approxsol )} \leq C_\kappa \, \eps_tÊ\,  \vert u \vert_{\ctbg  (\approxsol )}$ as promised. 
\end{proof}

\subsection{The Conclusion of the Proof}

The estimates for the proof of Theorem~\ref{th:1} are now all in place and the conclusion of the theorem becomes a simple verification of the conditions of the Inverse Function Theorem. We choose $\gamma \in (2-n,0)$. First, the linearization satisfies the estimate
\[
\vert \tilde{ \mathcal{L}}_ t \, u \vert_{\cobg(\approxsol)} \geq C_1 \,  \vert u \vert_{\ctbg(\approxsol)} \, ,
\]
by Proposition~\ref{thm:linest} where $C_1 >0$ is a constant independent of $t$ when $t$ is close enough to $t_*$.  Therefore by the Inverse Function Theorem of Section \ref{sec:deform} along with Proposition \ref{prop:nonlin}, a solution of the deformation problem can be found if 
\[
\big\vert \Phi_t (0) - H_t \big\vert_{\cobg(\approxsol)} \leq \tfrac{1}{2} C_1 \,  R
\]
where $R = \kappa \, \eps_t^{2-\gamma}$ and if 
\[
\big\vert\Dif \Phi_t (f) \, u - \tilde{\mathcal{L}}_t  \, u \big\vert_{\cobg (\approxsol)} \leq \tfrac{1}{2} \, C_1 \, \,  \vert u \vert_{\ctbg (\approxsol)} \, .
\]
But the second nonlinear estimate above shows that $ \big\vert\Dif \Phi_t (f) \, u - \tilde{\mathcal{L}}_t  \, u \big\vert_{\cobg (\approxsol)} \leq C_{\kappa} \, \eps_t \,  \vert u \vert_{\ctbg (\approxsol)}$ and Proposition \ref{prop:error} shows that $\big\vert\Phi_t (0) -  H_t \big\vert_{\cobg(\approxsol)} \leq   C_0 \, \eps_t^{2-\gamma}$.  Hence the two conditions above can always be met if $t$ is sufficiently close to $t_*$ and $\kappa$ is large enough to ensure $\kappa \, C_1  \geq 2 \, C_0$.  This concludes the proof of Theorem~\ref{th:1}. 

\hfill \qedsymbol

\bibliography{cmc}
\bibliographystyle{amsplain}

\end{document}